\documentclass[aop,citesort,MSNbibl,dvips]{arximspdf}
\usepackage{graphicx}

\doi{10.1214/10-AOP567}
\volume{39}
\issue{3}
\pubyear{2011}
\firstpage{1097}
\lastpage{1121}

\makeatletter

\newcommand{\OO}{{\O}}

\newcommand{\xrightarrow}[1]{\mathop{\hbox to 1cm{\rightarrowfill}}_{#1}}

\newtheorem{theorem}{Theorem}
\newtheorem{lemma}[theorem]{Lemma}
\newtheorem{proposition}[theorem]{Proposition}
\newtheorem{corollary}[theorem]{Corollary}

\newcommand{\EE}{{\mathbb E}}
\newcommand{\PP}{{\mathbb P}}

\newcommand{\Poisson}{\operatorname{Poisson}}
\newcommand{\Vect}{\operatorname{Vect}}
\newcommand{\rank}{\operatorname{rank}}

\newcommand{\ind}{{\mathbf{1}}}
\newcommand{\Dom}{\operatorname{Dom}}
\newcommand{\R}{{\mathbb R}}
\newcommand{\C}{{\mathbb C}}
\newcommand{\N}{{\mathbb N}}
\newcommand{\ibf}{{\mathbf i}}
\newcommand{\jbf}{{\mathbf j}}
\newcommand{\kbf}{{\mathbf k}}
\newcommand{\cG}{{\mathcal G}}
\newcommand{\cT}{{T}}
\newcommand{\cH}{{\mathcal H}}
\newcommand{\cP}{{\mathcal P}}
\newcommand{\cA}{{\mathcal A}}
\newcommand{\cB}{{\mathcal B}}
\newcommand{\tT}{\tilde{T}}
\newcommand{\cL}{{\mathcal L}}
\newcommand{\LR}{\operatorname{LR}}
\newcommand{\la}{\langle}
\newcommand{\ra}{\rangle}
\newcommand{\ox}{\overline{x}}

\makeatother

\begin{document}
\begin{frontmatter}

\title{The rank of diluted random graphs}
\runtitle{Rank of random graphs}

\begin{aug}
\author[A]{\fnms{Charles} \snm{Bordenave}\corref{}\ead
[label=e1]{charles.bordenave@math.univ-toulouse.fr}},
\author[B]{\fnms{Marc} \snm{Lelarge}\ead
[label=e2]{marc.lelarge@ens.fr}} and
\author[B]{\fnms{Justin} \snm{Salez}\ead[label=e3]{justin.salez@ens.fr}}
\runauthor{C. Bordenave, M. Lelarge and J. Salez}
\affiliation{Universit\'e Toulouse III, INRIA and INRIA}
\address[A]{C. Bordenave\\
CNRS UMR5219 and Institut\\
de Math\'ematiques de Toulouse\\
Universit\'e Toulouse III\\
France\\
\printead{e1}} %adresu isvedimo komanda gale!
\address[B]{M. Lelarge\\
J. Salez\\
D\'epartement d'Informatique, Projet TREC\\
INRIA-\'Ecole Normale Sup\'erieure\\
France\\
\printead{e2}\\
\phantom{E-mail: }\printead*{e3}}
\end{aug}

% HISTORY:
\received{\smonth{9} \syear{2009}}
\revised{\smonth{3} \syear{2010}}

% ABSTRACT
%
\begin{abstract}
We investigate the rank of the adjacency matrix of large diluted random
graphs: for a sequence of graphs $(G_n)_{n\geq0}$ converging locally
to a Galton--Watson tree $T$ (GWT), we provide an explicit formula for
the asymptotic multiplicity of the eigenvalue $0$ in terms of the
degree generating function $\varphi_*$ of $T$. In the first part, we
show that the adjacency operator associated with $T$ is always
self-adjoint; we analyze the associated spectral measure at the root
and characterize the distribution of its atomic mass at $0$. In the
second part, we establish a sufficient condition on $\varphi_*$ for the
expectation of this atomic mass to be precisely the normalized limit of
the dimension of the kernel of the adjacency matrices of $(G_n)_{n\geq
0}$. Our proofs borrow ideas from analysis of algorithms, functional
analysis, random matrix theory and statistical physics.
\end{abstract}

% KEYWORDS
%
\begin{keyword}[class=AMS]
\kwd[Primary ]{05C80}
\kwd{15A52}
\kwd[; secondary ]{47A10}.
\end{keyword}
\begin{keyword}
\kwd{Random graphs}
\kwd{adjacency matrix}
\kwd{random matrices}
\kwd{local weak convergence}
\kwd{Karp and Sipser algorithm}.
\end{keyword}

\end{frontmatter}

%s1 ###
\section{Introduction}\label{sec1}

In this paper we investigate asymptotical spectral properties of the adjacency
matrix of large random graphs. To motivate our work, let us briefly mention
its implications in the special case of Erd\H{o}s--R\'enyi random
graphs. Let $G_n = (V_n, E_n)$ be an Erd\H{o}s--R\'enyi graph with
connectivity $c>0$ on the vertex set $V_n = \{1, \ldots, n\}$. In
other words,
we let each pair of distinct vertices $ij$ belong to the edge-set $E_n$ with
probability $c / n$, independently of the other pairs. The adjacency
matrix $A_n$ of $G_n$ is the $n \times n$ symmetric matrix defined by
$(A_n)_{ij} = \ind((ij ) \in E_n)$. Let $\lambda_1 (A_n) \geq\cdots
\geq\lambda_n (A_n)$ denote the eigenvalues of $A_n$ (with
multiplicities) and
\[
\mu_n = \frac1 n \sum_{i=1}^n \delta_{\lambda_i (A_n)}
\]
denote the
spectral measure of $A_n$. Our main concern will be the rank of $A_n$
\[
\rank( A_n ) = n - \dim\operatorname{ker} (A_n) = n - n \mu_n ( \{
0 \}
).
\]
%
%With this notations in hands, our work will imply
%
\begin{theorem}\label{th1}
\textup{(i)}
There exists a deterministic symmetric measure $\mu$ such that, almost
surely, for the weak convergence of probability measures,
\[
\lim_{n\to\infty} \mu_n = \mu.
\]

\textup{(ii)}
Let $0 <q < 1$ be the smallest solution to $q = \exp( -c
\exp(-c q))$. Then almost surely,
\[
\lim_{n\to\infty} \mu_n (\{0\}) = \mu(\{0\}) = q + e^{-c q} + c q
e^{- c q}-1.
\]
In other words, almost surely,
%
%e1 ###
%
\begin{equation}
\label{eq:rk1}\lim_{n\to\infty}\frac{\rank(A_n)}{n} = 2-q - e^{-c
q} - c q e^{- c q}.
\end{equation}
\end{theorem}

Apart from an improvement of the convergence, part (i) is not new; the
convergence in probability was first rigorously proved by Khorunzhy,
Shcherbina and
Vengerovsky \cite{khorunzhy04} (for an alternative proof, see
\cite{resArXiv} [note that it only implies $\limsup_n \mu_n (\{0\})
\leq
\mu(\{0\})$]).

In the sparse case, that is, when the connectivity $c$ grows with $n$
like $a \log
n$, the rank of $A_n$ has been studied by Costello, Tao and Vu
\cite{ctv} and Costello and Vu \cite{costellovu}. Their results imply
that for $a >1$, with high probability $\dim\operatorname{ker}(A_n) =
0$ while for $0 < a <1$, $\dim\operatorname{ker}(A_n)$ is of order of
magnitude $n^{1-a}$. Our theorem answers one of their open questions in
\cite{costellovu}.

The formula (\ref{eq:rk1}) already appeared in a remarkable paper by
Karp and Sipser \cite{karpsipser} as the
asymptotic size of the number of vertices left unmatched by a maximum
matching of $G_n$.
To be more precise, the function $G\mapsto
\dim\operatorname{ker}(G)$ is easily checked to be invariant under
``leaf removal,'' that is, if $G'$ is the graph obtained from $G$ by
deleting a
leaf and its unique neighbor, then $\dim\operatorname{ker}(G')=\dim
\operatorname{ker}(G)$. Karp and Sipser
\cite{karpsipser} study the effect of iterating this leaf removal on
the random graph $G_n$ until only isolated vertices and a
``core'' with minimum degree at least $2$ remain.
They show that the asymptotic number of isolated vertices is
approximately $(2- q -e^{-c q} - c q e^{-c q})n$ as $n\to\infty$, and
that the size of the core is $o(n)$ when $c\leq e$. Thus,
(\ref{eq:rk1}) follows by additivity of $G\mapsto\dim\operatorname
{ker}(G)$ on
disjoint components, as observed by Bauer and Golinelli
\cite{bauerg01}.
However for $c>e$, the size of the core is not negligible and
the same argument only leads to the following inequality:
\[
\liminf_{n\to\infty}\frac{\dim\operatorname{ker}(A_n)}{n}\geq q +
e^{-c q} + c q e^{- c q} -1.
\]
Bauer and Golinelli \cite{bauerg01} conjecture that this lower bound
should be the actual limit for all $ c $, which is equivalent to saying
that asymptotically the
dimension of the kernel of the core is zero.
The proof of this conjecture follows from our work (see Section \ref
{rem:match}).
%We will also provide a new proof of Karp and Sipser Theorem.

Our results are not restricted to Erd\H{o}s--R\'enyi graphs. They will in
fact hold for any sequence $(G_n)_{n\geq1}$ of random graphs converging
locally to a rooted Galton--Watson tree (GWT), provided the latter satisfies
a certain degree condition. The precise definition of
local convergence is recalled in Section \ref{sec:LWC}. It was introduced
by Benjamini and Schramm \cite{bensch} and Aldous and Steele
\cite{aldste}. A rooted GWT (see~\cite{aldlyo}) is characterized by its
\textit{degree distribution} $F_*$, which can be any probability
measure with
finite mean on $\N$: the root $\OO$ has offspring distribution $F_*$
and all other genitors have offspring distribution
$F$, where for all $ k \geq1$, $F(k-1) = k F_* (k) / \sum_\ell\ell F_*
(\ell)$. In the case of Erd\H{o}s--R\'enyi graphs with connectivity
$c$, the
limiting tree is simply a GWT with degree distribution $F_*=\Poisson(c)$.

The adjacency operator $A$ of a GWT [$\cT=(V,E)$] is a densely defined
symmetric linear operator on the Hilbert space $\ell^2(V)$ defined for
$\ibf,\jbf$ in $V$ by
\[
\langle A e_\ibf, e_\jbf\rangle= \ind
(\ibf\jbf\in E),
\]
where for any $\ibf\in V$, $e_\ibf$ denotes the base function
$\jbf\in V\mapsto\ind(\jbf=\ibf)$.
As we will show, if $F_*$ has a finite second moment, then $A$ has almost
surely a unique self-adjoint extension, which we also denote by $A$.
Consequently, for any
unitary vector $\psi\in\Dom(A)$, the spectral theorem guarantees the
existence and uniqueness of a probability measure $\mu_{\psi}$ on $\R$,
called the spectral measure associated with $\psi$, such that for any
$k\geq0$,
\[
\langle A^k\psi,\psi\rangle= \int_\R x^k\,d\mu_{\psi}(x).
\]
In particular, we may consider the spectral measure $\mu_{T}$ associated
with the vector~$e_{\OO}$, where $\OO$ is the root of the rooted tree
$T$. Our first main result is an explicit formula for $\EE\mu_{T}(\{
0\})$, the
expected mass at zero of the spectral measure at the root $\OO$ of a
rooted GWT $\cT$.
\begin{theorem}\label{th:leaf}
Let $\cT$ be a GWT whose degree distribution $F_*$ has a finite second
moment, and let $\varphi_*$ be the generating function of $F_*$. Then,
$
\EE\mu_{T} (\{ 0 \}) = \max_{x\in[0,1]}M(x),
$
where
\[
M(x)=\varphi'_*(1)x\ox+\varphi_*(1-x)+\varphi_*(1-\ox)-1\qquad
\mbox{with } \ox=\varphi'_*(1-x)/\varphi'_*(1).
\]
\end{theorem}

\begin{figure}%[b]

\includegraphics{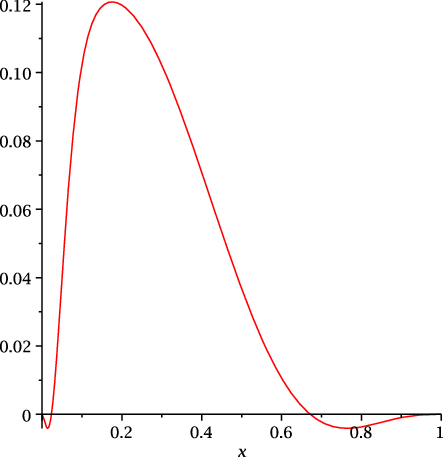}

\caption{Plot of $M$ for $\varphi_*(x)=\frac{d}{1+d}x^d
+\frac{1}{1+d}x^{d^3}$, with $d=3$.} \label{fig0}
\end{figure}

In the special case of regular trees, the measure $\mu_{T}$ can be
explicitely computed and
turns out to be absolutely continuous, so $\mu_{T}(\{0\})=0$. In
contrast, one
may construct GWTs with arbitrary large minimum degree and such that
$\EE\mu_T(\{0\})>0$. The following example is
taken form \cite{bofr} and is due to Picollelli and Molloy:
set $d\geq3$ and take $\varphi_*(x)=\frac{d}{1+d}x^d +
\frac{1}{1+d}x^{d^3}$. Figure \ref{fig0} gives\vspace*{1pt} a plot of $M$ for
the case $d=3$, showing that $\EE\mu_T(\{0\})>0$ in this case.

%
%f1 ###
%

When $F_*$ is a Poisson distribution with mean $c$, the corresponding
quantity $\max_{x\in[0,1]}M(x)$ is precisely (\ref{eq:rk1}), and it
already appeared in Zdeborov\'a
and M\'ezard~\cite{zdeborovamezard}, equation (38), as a ``cavity
method''
prediction for the limiting fraction of unmatched vertices in a maximum
matching.

To the best of our knowledge, the formula was unknown for general
GWTs. However, Bauer and Gollineli \cite{baugol00} have computed
explicitly the asymptotic rank of the uniform spanning tree on the complete
graph of size $n$. Also Bhamidi, Evans and Sen \cite{bhamidievanssen}
have recently analyzed the convergence of the spectrum of the adjacency
matrix of growing random trees.

Our second main result (Theorem \ref{th:main}) states that for any sequence
of random graphs $(G_n)_{n\geq0}$ converging locally in distribution
to a
GWT, we have $\lim_n n^{-1}\rank(A_n) = 1-\EE\mu_{T} (\{0\})$,
provided the first local extremum of the above function $x\mapsto
M(x)$ is a global maximum on $[0,1]$. We have left open the case where
the global maximum of $M$ is not the first local maximum (see Section
\ref{rem:match}).

Our detailed analysis of the atomic mass at $0$ of the limiting
spectral measure $\mu$ remains only a small achievement
for the global understanding of this measure. For example, for
Erd\H{o}s--R\'enyi graphs, the atomic part of $\mu$ is dense in $\R$,
and nothing is known on the mass of atoms apart $0$. There is also a
conjecture about the absolutely continuous part $\mu_{ac}$ of the measure
$\mu$: we say that $\mu$ has \textit{extended states} (resp., no extended
state) at $E \in\R$ if the partition function $x \mapsto\mu_{ac}
(-\infty, x)$ is differentiable at $x = E$ and its derivative is
positive (resp., null). This notion was introduced in mathematical physics
in the context of spectra of random Schr\"odinger operators; a recent
treatment can be found
in Aizenman, Sims and Warzel \cite{aizenman2006}. For
Erd\H{o}s--R\'enyi graphs, Bauer and Gollineli have conjectured that
$\mu$ has no extended state at $E = 0$ when $0< c \leq e$, and has extended
states at $E = 0$ when $c > e$. More generally, one may wonder whether
$\mu_{ac} = 0$ when $0< c \leq e$. Finally, the existence of a
singular continuous part in $\mu$ is apparently unknown.

The remainder of the paper is organized as follows: in Section
\ref{sec:GWTSA}, we analyze the adjacency operator of a GWT. In
Section \ref{sec:mu0}, we study $\mu_{T} (\{0\})$ and prove Theorem
\ref{th:leaf}. In Section \ref{sec:LWC}, we prove finally the
convergence of the spectrum of finite graphs and the convergence of
the rank. The proof of Theorem \ref{th1} is given in the
\hyperref[sec:ap]{Appendix}.
%We also include an appendix on random matchings in trees.

%s2 ###
\section{Locally finite graphs and their adjacency operators}
\label{sec:GWTSA}

A rooted graph is the pair formed by a graph $G$ with a distinguished
vertex $\OO\in V$, called the root. There is a canonical way to define
a distance on $V$: for each $u,v \in V$, the (graph)-distance is the
minimal length of a path from $u$ to $v$, if any, and $\infty$
otherwise. For a rooted graph $G$ with root $\OO$ and $t$ an integer,
we will denote by $(G)_t$ the rooted subgraph spanned by the vertices
at distance at most $t$ from the root. In all this section, we consider
a locally finite rooted graph $G=(V,E)$ with
root denoted by~$\OO$.
%If $\ibf\in V$, we use $|\ibf|$ for the distance from $\ibf$ to $\OO$.
%%, and we write $\jbf\succ\ibf$ to mean that $\ibf$ is an ancestor of $

%s2.1 ###
\subsection{Adjacency operator}
\label{sec:sa}

Consider the Hilbert space
\begin{eqnarray}
&&\ell^2(V)=\biggl\{\psi\dvtx V\to\C,
\sum_{\ibf\in V}|\psi(\ibf)|^2<\infty\biggr\}\nonumber\\
&&\eqntext{\mbox{with
inner product }\displaystyle\langle\psi,\phi\rangle=\sum_{\ibf\in
V}\psi(\ibf)\overline{\phi(\ibf)}.}
\end{eqnarray}
Denote by $H_0\subseteq\ell^2(V)$ the dense subspace of finitely
supported functions, and by $(e_\ibf)_{\ibf\in V}$ the canonical
orthonormal basis of $\ell^2(V)$, that is, $e_\ibf$ is the
coordinate function $\jbf\in V\mapsto\ind(\ibf=\jbf)$. By
definition, the adjacency operator $A$ of $G$ is the densely-defined linear
operator over $\ell^2(V)$ whose domain is $H_0$ and whose action on the
basis vector $e_\ibf,\ibf\in V$, is
\[
A e_\ibf= \sum_{\jbf\dvtx\ibf\jbf\in E}e_\jbf.
\]
Note that $A e_\ibf\in\ell^2(V)$ since $G$ is locally finite.
Moreover, for all $\ibf,\jbf\in V$,
\[
\la A e_\ibf, e_\jbf\ra= \ind\{\ibf\jbf\in E\} = \la A e_\jbf,
e_\ibf\ra.
\]
Therefore, the operator $A$ is symmetric, and we may now ask about the
self-adjointness of its closure, which is again denoted by $A$. The
answer of course depends upon $G$, but
here is a simple sufficient condition that should suit all our needs in
the present paper.

We define the boundary of a subset $ S\subseteq V$ as $\partial S=\{
\ibf\jbf\in E\dvtx\ibf\in S,\jbf\notin S\}$,
and the boundary degree $\Delta(\partial S)$ as the maximum number
of boundary edges that are adjacent to the same vertex.
\begin{proposition}
\label{pr:self-adj}
For $A$ to be self-adjoint, it is enough that $V$ admits an exhausting
sequence of finite subsets with bounded boundary degree:
\begin{enumerate}[(A)]
\item[(A)]\hypertarget{assumA}
There exist finite subsets $S_1,S_2,\ldots\subseteq V$
such that
\[
\bigcup_{n}S_n=V\quad\mbox{and}\quad\sup_{n}\Delta(\partial S_n)<\infty.
\]
\end{enumerate}
\end{proposition}
\begin{pf}
Denote by $A^*$ the adjoint of $A$. By the basic criterion for
self-adjointness (see, e.g., Reed and Simon \cite{reesim}, Theorem
VIII.3), it is enough to show that $0$ is the only vector
$\psi\in\Dom(A^*)$ satisfying $A^*\psi=\pm i\psi$. Consider such a
$\psi$
(let us treat, say, the $+i$ case), and define the following flow along
the oriented
edges of $G$:
\[
(\ibf\to\jbf) = \Im(\psi(\ibf)\overline{\psi(\jbf
)}) = -(\jbf\to\ibf),
\]
for all $\ibf\jbf\in E$. The amount of flow created at vertex $\ibf
\in V$
is then
\begin{eqnarray*}
\sum_{\jbf\dvtx\ibf\jbf\in E}(\ibf\to\jbf) &=& \Im
\biggl(\psi(\ibf)\sum_{\jbf\dvtx\ibf\jbf\in E}\overline{\psi
(\jbf)}\biggr)
=\Im\la A(\psi(\ibf)e_\ibf),\psi\ra\\
&=&\Im\la\psi(\ibf)e_\ibf,A^*\psi\ra
=|\psi(\ibf)|^2.
\end{eqnarray*}
Now, by anti-symmetry of the flow, the total amount of flow created inside
any finite subset $ S\subseteq V$ must equal the total amount of flow escaping
through the boundary $\partial S$
\[
\sum_{\ibf\in S}|\psi(\ibf)|^2=\sum_{\ibf\jbf\in\partial
S}(\ibf\to\jbf).
\]
Therefore, using $(\ibf\to\jbf)\leq|\psi(\ibf)||\psi(\jbf)|$ and
twice the Cauchy--Schwarz inequality, we find
\begin{eqnarray*}
\sum_{\ibf\in S}|\psi(\ibf)|^2 & \leq&
%& \leq&
\biggl(\sum_{\ibf\in\partial S^-}|\psi(\ibf)|^2\sum_{\ibf\in
\partial S^-}\biggl(\sum_{\jbf\in S^c\cap
N_\ibf}|\psi(\jbf)|\biggr)^2\biggr)^{1/2}\\
& \leq&
\Delta(\partial S)\biggl(\sum_{\ibf\in\partial S^-}|\psi(\ibf
)|^2\sum_{\jbf\in\partial S^+}|\psi(\jbf)|^2
\biggr)^{1/2},
\end{eqnarray*}
where we have written $N_\ibf$ for the set of neighbors of $\ibf$,
$\partial S^-$ and $\partial S^+$ for the sets of vertices $\partial
S\cap S$ and $\partial S\cap S^c$, respectively.
Finally, take
$ S= S_n$, and let $n\to\infty$: the exhaustivity
$\bigcup_n S_n=V$ ensures that the left-hand side tends to ${\sum_{\ibf
\in
V}}|\psi(\ibf)|^2=\Vert\psi\Vert^2$ and also that
\[
\sum_{\ibf\in\partial S_n^-}|\psi(\ibf)|^2\xrightarrow{n\to
\infty}^{}0
\quad\mbox{and}\quad\sum_{\jbf\in\partial S_n^+}|\psi(\jbf
)|^2\xrightarrow{n\to\infty}^{}0.
\]
Since $\sup_n\Delta(\partial S_n)<\infty$, the right-hand side
vanishes, and we obtain the desired $\Vert\psi\Vert=0$.
\end{pf}

%s2.2 ###
\subsection{Spectral measure}\label{sec:specm}

We now assume that the adjacency operator $A$ is
self-adjoint. The spectral theorem then guarantees the validity of the
\textit{Borel functional calculus} on $A$: any measurable function
$f\dvtx\R\to\C$ may now
be rigorously applied to the operator $A$ just as one would do with
polynomials. Denoting by $\mu_{G}$ the spectral measure
associated with the vector $e_{\OO}$, we may thus write
%
%e2 ###
%
\begin{equation}
\label{eq:mu0}\la f(A) e_{\OO},e_{\OO}\ra=\int_\R f(x)\,d\mu_{G}(x)
\end{equation}
for any $f\in\cL_\C(\mu_{G})$. Taking
$f(x)=x^n$ $(n\in\N)$, we obtain in particular
%
%e3 ###
%
\begin{eqnarray}
\label{eq:path}
\gamma_n&=&\la A^n e_{\OO},e_{\OO}\ra=\int x^n\,d\mu
_{G}(x) \nonumber\\[-8pt]\\[-8pt]
&=& \# \{\mbox{paths of
length $n$ from $\OO$ to $\OO$ in $G$}\}.\nonumber
\end{eqnarray}

Since $\Vert e_{\OO}\Vert=1$, the spectral measure $\mu_{T}$ is a probability
measure on $\R$. We will now study its Cauchy--Stieltjes
transform. By definition, the Cauchy--Stieltjes transform of a probability
measure $\mu$ on $\R$ is the holomorphic function $m_\mu$ defined on the
upper complex half-plane $\C_+$ by
\[
m_\mu\dvtx z \mapsto\int_\R\frac{d\mu(x) }{ x - z}.
\]
Note that $m_\mu$ belongs to the set $\cH$ of holomorphic
functions $f$ on $\C_+$ satisfying
\[
\forall z\in\C_+\qquad \Im f(z)\geq0\quad\mbox{and}\quad|f(z)| \leq
(\Im z)^{-1},
\]
which is compact in the normed space of holomorphic functions on $\C
^+$ (Montel's theorem).

Henceforth, we will assume that $G$ is a rooted tree $\cT$. We
write $\jbf\succ\ibf$ to mean that $\ibf\in V$ is an ancestor of
$\jbf\in
V$, and we let $\cT_\ibf$ be the subtree of $\cT$ restricted to $\{
\jbf\in V$, $\jbf\succeq\ibf\}$, rooted at $\ibf$. Its adjacency
operator $A_\ibf$ is the projection of $A$ on $\Vect(e_\jbf,
\jbf\succeq\ibf)$. Since it is also self-adjoint, we may
consider its spectral measure $\mu_{T_\ibf}$ associated with the
vector $e_{\ibf}$, and its Cauchy--Stieltjes transform
$m_{T_\ibf}$. The recursive structure of trees implies a simple
well-known recursion for the family $(m_{T_\ibf})_{\ibf\in V}$:
\begin{proposition}
The family $(m_{T_\ibf})_{\ibf\in V}$ is solution in $\cH^V$ to the
system of equations, for all $z\in\C_+$,
%
%e4 ###
%
\begin{equation}
\label{eq:recm} f _{\ibf}(z) =-\biggl(z+\sum_{\jbf\in D(\ibf)}
f_{\jbf}(z)\biggr)^{-1},
\end{equation}
where $D(\ibf)=\{\jbf\succ\ibf, |\jbf|=|\ibf|+1\}$ denotes the
set of
immediate children of $\ibf$.
\end{proposition}
\begin{pf}
As we will see, the recursion follows from a classical operator version
of the Schur
complement formula (see, e.g., Proposition 2.1 in Klein \cite{klein} for
a similar argument). We write the proof for completeness. Define the
operator $U$ on $\ell^2(V)$ by its matrix elements,
\[
\langle Ue_{\OO},e_\ibf\rangle=\langle Ue_\ibf,e_{\OO}\rangle
=1=\langle Ae_{\OO},e_{\ibf}\rangle,
\]
for all $\ibf\in D(\OO)$, and $\langle Ue_\jbf,e_\kbf\rangle=0$ otherwise.
We then have the following decomposition:
\[
A= U+ \bigoplus_{\ibf\in D(\OO)}A_\ibf,
\]
where $A_\ibf$, is the projection of $A$ on $V_\ibf=\Vect
(e_\jbf, \jbf\succeq\ibf)$.
Since $A$ and $\tilde{A}=\bigoplus_{\ibf\in D(\OO)}A_\ibf$ are
self-adjoint operators, their respective resolvents
\[
R\dvtx z\mapsto(A-zI)^{-1},\qquad \tilde{R}\dvtx z\mapsto
(\tilde{A}-zI)^{-1}
\]
are well defined on $\C_+$, and the resolvent
identity gives
%
%e5 ###
%
\begin{equation}
\label{eq:res}
R(z)U\tilde{R}(z) = R(z)-\tilde{R}(z).
\end{equation}
In particular, for all $\kbf\in V$,
\[
\langle R(z)U\tilde{R}(z)e_{\OO},e_\kbf\rangle= \langle R(z)
e_{\OO},e_\kbf\rangle-\langle\tilde{R}(z)e_{\OO},e_\kbf\rangle.
\]
Now, using the definition of $U$, we may expand the left-hand side as
\[
\biggl(\langle\tilde{R}(z) e_{\OO},e_\kbf\rangle\sum_{\ibf\in
D(\OO)}\langle R(z) e_{\OO},
e_{\ibf}\rangle\biggr) + \biggl(\langle R(z) e_{\OO},
e_{\OO}\rangle\sum_{\ibf\in D(\OO)} \langle\tilde{R}(z) e_{\ibf},
e_\kbf\rangle\biggr).
%&=& \langle R(z)
%e_{\OO},e_\kbf\rangle-\langle\tilde{R}(z)e_{\OO},e_\kbf\rangle.
\]
But $\tilde{R}(z)e_{\OO}=-z^{-1}e_{\OO}$ and each $V_\ibf, \ibf\in
D(\OO)$, is stable for $\tilde{R}$. Therefore, in the special case
where $\kbf= \OO$, the above equality simplifies into
\[
-\frac1z\sum_{\ibf\in D(\OO)}\langle R(z) e_{\OO},
e_{\ibf}\rangle= \langle R(z)
e_{\OO},e_{\OO}\rangle+ \frac1z,
\]
while for $\kbf\in D(\OO)$, it gives
\[
\langle R(z)e_{\OO},e_{\OO}\rangle\langle\tilde{R}(z)e_{\kbf},
e_{\kbf}\rangle= \langle R(z)e_{\OO},e_{\kbf}\rangle.
\]
Combining both, we finally obtain
\[
\langle R(z)e_{\OO},e_{\OO}\rangle= -\biggl(z+ \sum_{\ibf\in
D(\OO)}\la\tilde{R}(z)e_{\ibf},e_{\ibf}\ra\biggr)^{-1},
\]
which, by (\ref{eq:mu0}) with $f(x)=(x-z)^{-1}$, is precisely
\[
m_{T_{\OO}}(z) =-\biggl(z+ \sum_{\ibf\in D(\OO)}m_{T_\ibf}(z)\biggr)^{-1}.
\]
%
%Simply replace $\cT$ by $\cT_\ibf$.
\upqed\end{pf}

When $\cT$ is finite, the set of equations (\ref{eq:recm}) uniquely
determines the Cauchy--Stieltjes transforms $(m_{T_{\ibf}})_{\ibf\in
V}$, which can
be computed iteratively from the leaves up to the root. Under an extra
condition on $\cT$, this extends to the infinite case. Recall that
$(\cT)_n$ denote
the truncation of $\cT$ to the first $n$ generations.
In what follows, we will make the additional assumption
{\renewcommand{\theequation}{B}
\begin{equation}\label{equB}
\limsup_{n\to\infty}|\partial( \cT) _ n |^{
{1/n}}<\infty.
\end{equation}}
\begin{proposition}\label{prop:uni}
If $\cT$ satisfies assumption (\ref{equB}),
then $(m_{T_{\ibf}})_{\ibf\in V}$ is the unique solution in
$\cH^V$ to the system of equations (\ref{eq:recm}), and for all $\ibf
\in
V$,
%
%e6 ###
%
\setcounter{equation}{5}
\begin{equation}
\label{eq: convmeasure}
m_{T_{\ibf}} = \lim_{n\to\infty}m_{(T_\ibf)_n},
\end{equation}
in the sense of compact convergence on $\C_+$.
\end{proposition}
\begin{pf}
If $(f_{\ibf})_{\ibf\in V}\in\cH^V$ and $
(g_{\ibf})_{\ibf\in
V}\in\cH^V$ are solutions to the system of equations (\ref
{eq:recm}), then
we can write, for all $\ibf\in V$, $z\in\C_+$,
\begin{eqnarray*}
|f_\ibf(z)-g_\ibf(z)| &=& \biggl|\frac{\sum_{\jbf\in D(\ibf
)}(f_{\jbf}(z)-g_\jbf(z))}{(z+\sum_{\jbf\in D(\ibf)}f_\jbf
(z))(z+\sum_{\jbf\in
D(\ibf)}g_\jbf(z))}\biggr|\\
&\leq&{\frac{1}{(\Im(z))^{2}} \sum
_{\jbf\in
D(\ibf)}}|f_\jbf(z)-g_\jbf(z)|.
\end{eqnarray*}
Iterating this $n$ times, and then using the uniform bound $|f_\jbf
(z)-g_\jbf(z)|\leq2\times(\Im(z))^{-1}$, we obtain
\[
|f_\ibf(z)-g_\ibf(z)|\leq{\frac{1}{(\Im(z))^{2n}} \sum_{\jbf\in
\partial( \cT_\ibf)_ n }}|f_\jbf(z)-g_\jbf(z)|\leq\frac{2 |
\partial( \cT_\ibf)_ n |}{(\Im(z))^{2n+1}}.
\]
Therefore, we see that under assumption (\ref{equB}),
\[
\forall\ibf\in
V\qquad |f_\ibf(z)-g_\ibf(z)|=0
\]
as soon as $\Im(z)$ is sufficiently large,
hence for all $z\in\C_+$ by holomorphy. Finally, denote by $M_n$ the
denumerable vector of holomorphic functions $(m_{(T_\ibf)_n})_{\ibf
\in
V}\in\cH^V$. Since $\cH$ is compact, the sequence $(M_n)_{n\geq0}$ is
relatively compact, and since each vector $M_n$ satisfies the partial
set of equations (\ref{eq:recm}) corresponding
to the truncated tree $(T)_n$, any limit point $M_\infty$ must satisfy the
global set of equations (\ref{eq:recm}) corresponding to the full tree
$\cT$, so $M_\infty$ is nothing but $(m_{T_{\ibf}})_{\ibf\in V}$.
Therefore, the
sequence of vectors $(M_n)_{n\geq0}$ converges to $M_\infty
=(m_{T_{\ibf}})_{\ibf\in V}$,
and this is exactly (\ref{eq: convmeasure}).
\end{pf}

%s2.3 ###
\subsection{Atomic mass at zero}
\label{sec:mu0}

Our goal here is to characterize $\mu_{T}(\{0\})$, the atomic mass at
zero of the
spectral measure $\mu_{T}$.
\begin{proposition}\label{lem:monot}
If $\cT$ satisfies assumption (\ref{equB}),
then
the family $(\mu_{T_\ibf} (\{0\}))_{\ibf\in V}$ is the
largest solution in $[0,1]^V$ to the system of equations
%
%e7 ###
%
\begin{equation}
\label{eq:recmu0}x_\ibf=\biggl( 1+\sum_{\jbf\in
D(\ibf)}\biggl(\sum_{\kbf\in D(\jbf)} x_\kbf\biggr)^{-1}\biggr)^{-1},
\end{equation}
with the conventions $1/0 = \infty$ and $1/\infty= 0$.
\end{proposition}
\begin{pf}
Since $\cT$ is acyclic, (\ref{eq:path}) ensures that the measures
$\mu_{T_\ibf}, \ibf\in V$, are symmetric. Therefore, for all $t>0,
\ibf\in V$
\[
m_{T_{\ibf}}(it) = \int_\R\frac{x}{x^2+t^2}\,d\mu_{T_\ibf}
(x)+i\int_\R\frac{t}{x^2+t^2} \,d\mu_{T_\ibf} (x)= i\int_\R\frac
{t}{x^2+t^2} \,d\mu_{T_\ibf} (x).
\]
Hence, if we define $h_{T_{\ibf}}(t):=-it m_{T_\ibf}(it)\in[0,1]$,
then by the dominated convergence theorem,
\[
h_{T_{\ibf}}(t) = \int_\R\frac{t^2\,d\mu_{T_\ibf}(x)}{x^2+t^2}
\xrightarrow{t\to0}^{} \mu_{T_\ibf} (\{0\}).
\]
But, iterating once equation (\ref{eq:recm}), we get
%
%e8 ###
%
\begin{equation}
\label{eq:recth}h_{T_{\ibf}}(t) =\biggl( 1+\sum_{\jbf\in
D(\ibf)}\biggl(t^2+\sum_{\kbf\in D(\jbf)} h_{T_\kbf} (t)
\biggr)^{-1}\biggr)^{-1},
\end{equation}
so that letting $t\to0$ yields exactly that $(\mu_{T_\ibf} (\{0\}
))_{\ibf\in
V}$ must satisfy (\ref{eq:recmu0}).
%the following equation (by continuity of $x \mapsto x^{-1}$ on $[0,

Again, when the rooted tree $\cT$ is finite, this recursion
characterizes the
family $(\mu_{T_\ibf} (\{0\}))_{\ibf\in V}$, since it
can be computed
iteratively from the leaves up to the root. However, when $\cT$ is
infinite, (\ref{eq:recmu0}) may admit several other solutions.
Fortunately, among all of them, $(\mu_{T_\ibf} (\{0\})
)_{\ibf\in V}$
is always the largest. To see why,\vspace*{2pt} consider any solution $(x_\ibf
)_{\ibf\in
\cT}\in[0,1]^V$. Fixing $t>0$, let us show by induction that for all
$n\in\N$,
%
%e9 ###
%
\begin{equation}
\label{eq:recur}
\forall\ibf\in V\qquad x_\ibf\leq h_{(T_\ibf)_{2n} }(t):=-it
m_{(T_\ibf)_{2n} }(it).
\end{equation}
This will conclude our proof since we may then let $n\to\infty$ to obtain
$x_\ibf\leq h_{T_{\ibf}}(t)$ by Proposition \ref{prop:uni}, and let
finally $t\to0$ to reach the desired $x_\ibf\leq\mu_{T_\ibf} (\{0\}
)$. The
base case $n=0$ is trivial because the right-hand equals $1$. Now, if
(\ref{eq:recur}) holds for some $n\in\N$, then for all $\ibf\in
V$,
\begin{eqnarray*}
x_\ibf &=& \biggl( 1+\sum_{\jbf\in
D(\ibf)}\biggl(\sum_{\kbf\in D(\jbf)} x_\kbf\biggr)^{-1}
\biggr)^{-1} \\
&\leq& \biggl( 1+\sum_{\jbf\in
D(\ibf)}\biggl(t^2+\sum_{\kbf\in D(\jbf)} h_{(T_\kbf)_{2n}
}(t)\biggr)^{-1}\biggr)^{-1}
=h_{(T_\ibf)_{2n+2} }(t),
\end{eqnarray*}
where the first equality follows from the fact that $(x_\ibf)_{\ibf
\in
\cT}$ satisfies (\ref{eq:recmu0}), the middle inequality from the
induction hypothesis, and the last equality from (\ref{eq:recth})
applied to $(\cT_\ibf)_{2n+2}$.
%
%The second point of the proposition follows from the fact that
%$\mu_{T_\ibf}(\{0\})\in[0,1]$ and the monotonicity of (
\end{pf}

%s2.4 ###
\subsection{Galton--Watson trees}\label{sec:GWT}
We now apply the above results to Galton--Watson trees.
Let $F_*$ be a distribution on $\N$ with finite mean, and let $\cT$
be a GWT with \textit{degree} distribution $F_*$, that is, a random
locally finite rooted
tree obtained by a Galton--Watson branching process where the root has
offspring distribution $F_*$, and all other genitors have
offspring distribution $F$, where
%
%e10 ###
%
\begin{equation}\label{eq:F}
\forall k\geq1\qquad F(k-1) = k F_* (k) \big/ \sum_\ell\ell F_* (\ell).
\end{equation}
In the rest of this paper, we will make the following second moment
assumption on the distribution $F_*\dvtx\sum_k k^2F_*(k) < \infty$, or
equivalently $\sum_k k F(k) < \infty$.
It is in fact a sufficient condition for all the previous results to hold
almost surely.
\begin{proposition}\label{prop:sa}
If $F_*$ has a finite second moment, then $\cT$ satisfies
\hyperlink{assumA}{\textup{(A)}} and~(\ref{equB}) with probability one.
In particular, the adjacency operator $A$ is almost surely
self-adjoint, and the atomic mass at zero of the spectral measure at
the root of $\cT$ is characterized by the fixed-point equation
(\ref{eq:recmu0}).
\end{proposition}
\begin{pf}
Let $N$ denote a generic random variable with law $F$. For (\ref{equB}), it is
well known (and easy to check by a martingale argument) that the size
of the $n$th
generation in a GWT with offspring distribution $F$ behaves like $\EE
^n N$ as $n\to\infty$, in the
precise sense that almost surely,
${ n^{-1}\log}|\partial(T)_n|\to\EE N$, which is finite by assumption.
As far as \hyperlink{assumA}{(A)} is concerned now, if $\cT$ is finite there is nothing to
do. Now if $\cT$ is infinite, we build an exhausting sequence of finite
vertex subsets with uniformly bounded boundary degree as follows:
the finite first moment assumption on $F$ guarantees the existence of
a large enough integer $\kappa\geq1$ so that
%
%e11 ###
%
\begin{equation}
\label{eq:sumtrunc}
\sum_{k\geq\kappa}kF(k)<1.
\end{equation}
For each vertex of $T$, color it in red if it has less than $\kappa$
children and in blue otherwise. If the root\vadjust{\goodbreak} $\OO$ is red, set
$ S _1=\{\OO\}$. Otherwise, the connected blue component containing the root
is a GWT with average offspring $\sum_{k\geq\kappa}kF(k)<1$, so it is
almost-surely finite, and we define $ S _1$ as the set of its vertices,
together with their (red) external boundary vertices. Now for each
external boundary vertex $\ibf\in\partial S _1^+$, we repeat the procedure
on the subtree $\cT_\ibf$, and we define $ S _2$ as the union of $ S _1$
and all the resulting subsets. Iterating this procedure, we obtain an
exhaustive sequence of subsets $ S _1, S _2,\ldots\subseteq V$ whose
boundary degree satisfies by construction $\Delta(\partial S
_n)=\kappa$,
which is exactly \hyperlink{assumA}{(A)}.
\end{pf}

Owing to the recursive distributional nature of GWTs, the set of equations
(\ref{eq:recmu0}) defining $\mu_{T}(\{0\})$ takes the much nicer form
of a
Recursive distributional equation (RDE), which we now make explicit.
We denote $\cP(\N)$ (resp., $\cP([0,1])$) the space of probability
distributions on $\N$ ($[0,1]$, resp.).
Given $F,{F'}\in\cP(\N)$ and $\nu\in\cP([0,1])$, we denote
by $\Theta_{F,{F'}}(\nu)$ the distribution of the $[0,1]$-valued r.v.
%
%e12 ###
%
\begin{equation}
\label{eq:operator}
Y=\frac{1}{1+\sum_{i=1}^{N}(\sum_{j=1}^{{N_i'}}X_{ij})^{-1}},
\end{equation}
where ${N}\sim F$, ${N_i'}\sim{F'}$ and ${X}_{ij}\sim\nu$, all of
them being independent. With this notation in hand, the previous
result implies the following: if $F^*$ has a finite second moment,
then $\mu_{T}(\{0\})$ has distribution $\Theta_{F_*,F}(\nu_0^*)$,
where $F$ is given by (\ref{eq:F}) and $\nu^*_0$ is the largest
solution to the RDE
%
%e13 ###
%
\begin{equation}
\label{eq:rdezero}
\nu^*_0=\Theta_{F,F}(\nu^*_0).
\end{equation}

The remainder of this section is dedicated to solving
(\ref{eq:rdezero}) when $F_*$ has a finite second moment. We will
assume that $F_*(0)+F_*(1)< 1$; otherwise $F=\delta_0$ and
$\nu_0^*=\delta_1$ is clearly the only solution to (\ref{eq:rdezero}).
We let $\varphi_*(z)=\sum_{n\geq0}F_*(n)z^n$ be the generating
function of $F_*$. For any $x\in[0,1]$, we set $\overline{x}=\varphi
'_*(1-x)/\varphi'_*(1)$, and we define
\[
M(x)=\varphi'_*(1)x\overline{x}+\varphi_*(1-x)+\varphi
_*(1-\overline{x})-1.
\]
Observe that $M'(x)=\varphi''_*(1-x)(\overline{\overline x} -x
)$, and therefore any $x \in[0,1]$ where $M$ admits a local
extremum must satisfy $x=\overline{\overline{x}}$. We will say that
$M$ admits a historical record at $x$ if $x=\overline{\overline{x}}$
and $M(x)>M(y)$ for any $0\leq y < x$. Since $[0,1]$ is compact and $M$
is analytic, there are only finitely many such records. In fact, they
are in one-to-one correspondence with the solutions to the RDE (\ref
{eq:rdezero}).
\begin{theorem}
\label{th:rde}
If $p_1<\cdots<p_r$ are the locations of the historical records of
$M$, then the RDE (\ref{eq:rdezero}) admits exactly $r$ solutions;
moreover, these solutions can be stochastically ordered, say $\nu
_{1}<\cdots<\nu_{r}$, and for any $i\in\{1,\ldots,r\}$:
\begin{enumerate}[(ii)]
\item[(i)] $\nu_{i}(\{0\}^c)=p_i$;
\item[(ii)] $\Theta_{F_*,F}(\nu_{i})$ has mean $M(p_i)$.
\end{enumerate}
In particular, $\EE[\mu_{T}(\{0\})]=\max_{x\in[0,1]}M(x)$.
\end{theorem}

It now remains to prove Theorem \ref{th:rde}. The space $\cP
([0,1])$ is naturally equipped with:
\begin{itemize}
\item[-] a natural topology, which is that of weak convergence,
\begin{eqnarray}
&&\mu_n\xrightarrow{n\to\infty}^{}\mu\quad\Longleftrightarrow\quad\int
_{}\varphi \,d\mu_n\xrightarrow{n\to\infty}^{}\int_{}\varphi
\,d\mu\nonumber\\
&&\eqntext{\mbox{for any continuous function }\varphi\dvtx[0,1]\to\R;}
\end{eqnarray}
\item[-] a natural order, which is that of stochastic domination,
\begin{eqnarray}
&&\mu_1\leq\mu_2 \quad\Longleftrightarrow\quad\int_{}\varphi \,d\mu_1\leq\int
_{}\varphi \,d\mu_2\nonumber\\
&&\eqntext{\mbox{for any continuous, increasing function
}\varphi\dvtx[0,1]\to\R.}
\end{eqnarray}
\end{itemize}
The proof is based on two lemmas, the first one being straightforward.
\begin{lemma}
\label{lm:contincr}
For any $F,F'\in\cP(\N)\setminus\{\delta_0\}$, $\Theta_{F,F'}$ is
continuous and strictly increasing on $\cP([0,1])$.
\end{lemma}
%
%It follows directly from the fact that, for any $n\geq0$ and any $n_1,
%$$x\mapsto\frac{1}{1+\sum_{i=1}^n(\sum_{j=1}^{n_i}x_{ij}
%)^{-1}}$$ is continuous and increasing from $[0,1]^{n_1+
%
\begin{lemma}
\label{lm:pF(p)}
For any $\nu\in\cP([0,1])$, letting $p=\nu(\{0\}
^c)$, we have:
\begin{enumerate}[(iii)]
\item[(i)] $\Theta_{F,F}(\nu)(\{0\}^c)=\overline
{\overline{p}}$;
\item[(ii)] if $\Theta_{F,F}(\nu)\leq\nu$, then the mean of
$\Theta_{F_*,F}(\nu)$ is at least $M(p)$;
\item[(iii)] if $\Theta_{F,F}(\nu)\geq\nu$, then the mean of
$\Theta_{F_*,F}(\nu)$ is at most $M(p)$.
\end{enumerate}
In particular, if $\nu$ is a fixed point of $\Theta_{F,F}$, then
$p=\overline{\overline{p}}$ and $\Theta_{F_*,F}(\nu)$ has mean $M(p)$.
\end{lemma}
\begin{pf}
In (\ref{eq:operator}) it is clear that $Y>0$ if and only if
for any $i\in\{1,\ldots,N\}$, there exists $j\in\{1,\ldots,{N_i'}\}
$ such that $X_{ij}>0$. Denoting by $\varphi$ the generating function
of~$F$, this rewrites
\[
\Theta_{F,F}(\nu)(\{0\}^c)=\varphi\bigl(1-{\varphi
}\bigl(1-\nu(\{0\}^c)\bigr)\bigr).
\]
But from (\ref{eq:F}) it follows that ${\varphi}(\cdot)=\varphi
'_*(\cdot)/\varphi'_*(1)$, that is, $\varphi(1-x
)=\overline x$, hence the first result.

Now let $X\sim\nu$, $N_*\sim F_*$, $N\sim F$, and let $S,S_1,\ldots$
have the distribution of the sum of $N$ i.i.d. copies of $X$, all these
variables being independent. Then, $\Theta_{F_*,F}(\nu)$ has mean
\begin{eqnarray*}
\EE\biggl[\frac{1}{1+\sum_{i=1}^{N_*} S_i^{-1}}\biggr]
& = &\EE\biggl[\biggl(1-\frac{\sum_{i=1}^{N_*} S_i^{-1}}{1+\sum
_{i=1}^{N_*} S_i^{-1}}\biggr)\mathbf1_{\{\forall i=1,\ldots,
{N_*}, S_i>0\}}\biggr]\\
& = &
\varphi_*(1-\overline{p})\\
&&{} -\varphi'_*(1)
\EE\biggl[\frac
{S^{-1}}{S^{-1}+1+\sum_{i=1}^{ {N}}S_i^{-1}}\mathbf1_{\{
S>0,\forall i=1,\ldots,\widehat{N}_*, S_i>0\}}\biggr]\\
& = &
\varphi_*(1-\overline{p})-\varphi'_*(1)\EE\biggl[\frac
{Y}{Y+S}\mathbf1_{\{S>0\}}\biggr],
\end{eqnarray*}
where the second and last lines follow from (\ref{eq:F}) and $Y\sim
\Theta_{F,F}(\nu)$, respectively. Now, for any $s>0$, $x\mapsto\frac
{x}{x+s}$ is increasing, and hence, depending on whether $\Theta
_{F,F}(\nu)\geq\nu$ or $\Theta_{F,F}(\nu)\leq\nu$, $\Theta
_{F_*,F}(\nu)$ has mean at most/least
%
%e14 ###
%
\begin{eqnarray}
\label{eq:calc}\qquad
&&\varphi_*(1-\overline{p})-\varphi'_*(1)\EE\biggl[\frac
{X}{X+S}\mathbf1_{\{S>0\}}\biggr]\nonumber\\[-8pt]\\[-8pt]
&&\qquad=
\varphi_*(1-\overline{p})-p\varphi'_*(1)\EE\biggl[\frac
{1}{1+\widehat N}\mathbf1_{\{\widehat N\geq1\}}\biggr]
\qquad\mbox{with }
\widehat N=\sum_{i=1}^{N}\mathbf1_{\{X_i>0\}}.\nonumber
\end{eqnarray}
%
%where we used exchangeability as follows:
%&&\EE[\frac{X}{X+\sum_{i=1}^{N}X_i}\mathbf1_{\{
%&=& \EE[ \mathbf1_{\{\widehat N=\sum_{i=1}^{N}\mathbf1_{\{X_i>0
%[ \frac{X_k}{\sum_{\ell=1}^{\widehat N+1}X_\ell}| \widehat N, X_1,
%&=& \EE[ \mathbf1_{\{\widehat N\geq1\}} \frac{1}{\widehat N +1}
%N+1}X_\ell}| \widehat N, X_1,\dots X_{\widehat N+1}>0]\PP(X_k>0)
%]\\
%&=& p\EE[ \mathbf1_{\{\widehat N\geq1\}} \frac{1}{\widehat N +1}
%]

But using the definition (\ref{eq:F}) and the well-known identity
$(n+1){n\choose d}=(d+1){n+1\choose d+1}$, one can easily check that
\begin{eqnarray*}
&&\varphi_*(1-\overline{p})-p\varphi'_*(1)\EE\biggl[\frac
{1}{1+\widehat N}\mathbf1_{\{\widehat N\geq1\}}\biggr]\\
&&\qquad= \varphi_*(1-\overline{p})- p\varphi'_*(1)\sum_{n\geq1}F(n)\sum
_{d=1}^n\pmatrix{n\cr d}\frac{p^d(1-p)^{n-d}}{d+1}
%& = & \sum_{n\geq1}\pi_{n+1}\sum_{d=1}^{n+1}{n+1\choose
%d+1}p^{d+1}(1-p)^{n-d}\\
%& = & \sum_{n\geq2}\pi_{n}(1-(1-p)^{n}-np(1-p)^{n-1})\\
\\
&&\qquad= M(p).
\end{eqnarray*}
\upqed\end{pf}

We now have all the ingredients we need to prove Theorem \ref{th:rde}.
\begin{pf*}{Proof of Theorem \ref{th:rde}}
Let $p\in[0,1]$ such that $\overline{\overline{p}}=p$, and define
$\nu_0=\operatorname{Bernoulli}(p)$.
%and then by induction $\mu_{k+1}=\Theta_{\widehat\pi,\widehat\pi}(
From Lemma \ref{lm:pF(p)} we know that $\Theta_{F,F}(\nu
_0)(\{0\}^c)=p$, and since $\operatorname{Bernoulli}(p)$
is the largest element of $\cP([0,1])$ putting mass $p$ on $\{0\}^c$,
we have $\Theta_{F,F}(\nu_0)\leq\nu_0$. Immediately,
Lemma \ref{lm:contincr} guarantees that the limit
\[
\nu_\infty=\lim_{k\to\infty}\searrow\Theta^k_{F,F}(\nu
_0)
\]
exists in $\cP([0,1])$ and is a fixed point of $\Theta
_{F,F}$. Moreover, by Fatou's lemma, the number $p_\infty=\nu_\infty
(\{0\}^c)$ must satisfy $p_\infty\leq p$. But then the
mean of $\Theta_{F_*,F}(\nu_\infty)$ must be both:
\begin{itemize}
\item[-] equal to $M(p_\infty)$ by Lemma \ref{lm:pF(p)} with $\nu
=\nu_\infty$ and
\item[-] at least $M(p)$ since $\forall k\geq0$, the mean of $\Theta
_{F_*,F}(\Theta^k_{F,F}(\mu_0))$ is at least $M(p)$
[Lemma \ref{lm:pF(p)} with $\nu=\Theta^k_{F,F}(\nu_0)$].
\end{itemize}
We have just shown both $M(p)\leq M(p_\infty)$ and $p_\infty\leq p$.
From this, we will now deduce the one-to-one correspondence between
historical records of $M$ and fixed points of $\Theta_{F,F}$. We treat
each inclusion separately:
\begin{itemize}
\item[-] If $M$ admits a historical record at $p$, then clearly
$p_\infty=p$, so $\nu_\infty$ is a fixed point satisfying $\nu
_\infty(\{0\}^c)=p$.% and $\Theta_{\pi,\widehat\pi}(\mu_
\item[-] Conversely, considering a fixed point $\nu$ with $\nu
(\{0\}^c)=p$, we want to deduce that $M$ admits a historical
record at $p$. We first claim that $\nu$ is the above defined limit
$\nu_\infty$. Indeed, $\nu\leq\operatorname{Bernoulli}(p)$ implies $\nu
\leq\nu_\infty$ ($\Theta_{F,F}$ is increasing), and in particular
$p\leq p_\infty$. Therefore, $p=p_\infty$ and $M(p)=M(p_\infty)$. In
other words, the two ordered distributions $\Theta_{F_*,F}(\nu)\leq
\Theta_{F_*,F}(\nu_\infty)$ share the same mean and hence are equal.
This ensures $\nu=\nu_\infty$. Now, if $q<p$ is any historical
record location, we know from part 1 that
\[
\lambda_\infty=\lim_{k\to\infty}\searrow\Theta_{F,F}^{k}
(\operatorname{Bernoulli}(q))
\]
is a fixed point of $\Theta_{F,F}$ satisfying $\lambda_\infty
(\{0\}^c)=q$. But $q<p$, so $\operatorname{Bernoulli}(q)<
\operatorname{Bernoulli}(p)$, hence $\lambda_\infty\leq\nu_\infty$. Moreover,
this limit inequality is strict because $\lambda_\infty(\{0\}
^c) = q < p = \nu_\infty(\{0\}^c)$. Consequently,
$\Theta_{F_*,F}(\lambda_\infty)< \Theta_{F_*,F}(\nu_\infty)$ and
taking expectations, $M(q)< M(p)$. Thus, $M$ admits a historical
record at $p$.\qed
\end{itemize}
\noqed\end{pf*}

%s3 ###
\section{Convergence of the spectral measure}
\label{sec:LWC}

%s3.1 ###
\subsection{Local convergence of rooted graphs}

In this paragraph, we briefly recall the framework of local convergence
introduced by Benjamini and Schramm \cite{bensch} and Aldous
and Steele \cite{aldste} (see also Aldous and Lyons \cite{aldlyo}).

% Let $G = (V,E)$ be a simple, locally finite graph. A path of length
%$n$ from $u$ to $v$ is a sequence of $n+1$ vertices $u_0,u_1,
%$k=0\dots n-1$. There is a canonical way to define a distance $d_G$ on
%$V$: for each $u,v \in V$, $d_G (u , v)$ is the minimal length of a
%path from $u$ to $v$, if any, and $\infty$ otherwise.
%The ball of radius $t$ centered at $u$ is then naturally $$B_G (u, t)
%:= \{ v \in V: d_G (u , v) < t \}.$$

%We recall that a rooted graph $(G, \OO)$ is the pair formed by a graph
%$G$ with a distinguished vertex $\OO\in V$, called the root.
We recall that for integer $t$, $(G)_t$ is the rooted subgraph spanned
by the vertices at distance at most $t$ from the root.
%For any $t > 0$, we let $[G,\OO][t]$ denote the graph whose vertex-set
%is $B_G (\OO, t)$ and whose edges are those of $G$ that have both
%end-points in $B_G (\OO, t)$.
We consider the set $\cG_*$ of all locally finite, connected rooted
graphs, taken up to root-preserving isomorphism. With the terminology
of combinatorics, $\cG_*$ is the set of rooted unlabeled connected
locally finite graphs. We define a metric on $\cG_*$ by letting the
distance between
two rooted graphs $G_1$ and $G_2 $ be $1 / (1 + T)$, where $T$ is the supremum
of those $t\geq0$ such that there exists a root-preserving isomorphism
from $(G_1)_t$ to $(G_2)_t$.
Convergence with respect to this metric is called \textit{local} convergence.

This makes $\cG_*$ into a separable and complete metric space (see
Section 2 in \cite{aldlyo}).
In particular, we can endow $\cG_*$ with its Borel $\sigma$-algebra and
speak about weak convergence of random elements in $\cG_*$.
Specifically, a sequence of
probability distributions $\rho_1,\rho_2,\ldots$ on $\cG_*$
converges weakly to a probability distribution $\rho$, denoted by
$\rho_n \Longrightarrow\rho$, if
\[
\int_{\cG^*}f\,d\rho_n\xrightarrow{n\to\infty}^{} \int_{\cG
^*}f\,d\rho
\]
for all bounded continuous function $f\dvtx\cG_*\to\R$.
This is called the \textit{local weak} convergence.

Let us finally mention three important examples of random graph
sequences that
converge locally weakly to Galton--Watson trees. The
Erd\H{o}s--R\'enyi graphs with connectivity $c$ on the vertex set
$\{1,\ldots,n\}$, rooted at $\OO=1$ converges locally weakly to the
GWT with degree distribution
$\Poisson(c)$. The uniform $k$-regular ($k\geq2$) graph on
$\{1,\ldots,n\}$, rooted at $\OO=1$, converges weakly to the infinite
$k$-regular tree. More generally, if $F_*$ is a degree distribution on
$\N$
with finite mean, the random graph-sequence with asymptotic degree
distribution $F_*$ converges to the GWT with degree distribution
$F_*$. Note that in the above examples, the vertices are exchangeable and
the choice $\OO=1$ is arbitrary: equivalently, we could have chosen $\OO
$ uniformly at
random among all vertices, independently of the edge structure.

%s3.2 ###
\subsection{Continuity of the spectral measure}

Since the elements of $\cG^*$ have countably many vertices and are only
considered up to isomorphism, we may without loss of generalities embed all
vertices into the same, fixed generic vertex set $V$, say the set of finite
words over integers: the root is represented by the empty-word~$\OO$, and
vertices at distance $t$ from the root are represented by word of length
$t$ in the usual way. All adjacency operators can thus be viewed as
acting on the same Hilbert space
$\ell^2(V)$, their action being defined as zero on the orthogonal
complement of the subspace spanned by their vertices. Note that this
does not affect
the spectral measure at the root $\mu_{T}$.

If $(G_n)$ is a converging sequence in $\cG_*$, say to $G\in\cG_*$,
we may
even relabel the vertices in a consistent way so that the
root-preserving isomorphisms
appearing in the definition of local convergence become identities:
for every $t\in\N$, there exists $n_t\in\N$ such that
%
%e15 ###
%
\begin{equation}
\label{eq:cvloc}
n\geq n_t \quad\Longrightarrow\quad(G_n)_t=(G)_t.
\end{equation}

Fixing a word $\ibf\in V$, and setting $t$ equal $1$ plus the distance
from $\ibf$ to the root above, we obtain that for all $n\geq n_t$,
$\ibf$ is a vertex of $G_n$ if and only if it is a vertex of $G$, and
in that case its neighbors in $G_n$ are exactly its neighbors in
$G$. In other words, $A_ne_\ibf= Ae_\ibf$. By linearity, it follows
that any finitely supported vector $\psi\dvtx V\to\C$ must satisfy
\[
A_n\psi\xrightarrow{n\to\infty}^{\ell^2(V)} A\psi,
\]
and since those $\psi$ are dense in ${\ell^2(V)}$, Theorem
VIII.25(a) in Reed and Simon \cite{reesim} guarantees that $A_n\to A$
in the
strong resolvent sense, provided of course that $A,A_1,\ldots$ are
self-adjoint. In particular, this implies the weak convergence of the
corresponding spectral measures at the root and the compact convergence of
their associated Cauchy--Stieltjes transforms,
\[
m_{G_n}\xrightarrow{n\to\infty}^{\cH}m_{G} \quad\mbox{and}\quad \mu
_{G_n}\xrightarrow{n\to\infty}^{\cP(\R)}\mu_{G}.
\]
Note that this last statement does not depend anymore on the way
$G,G_1,\ldots$ have been embedded. We have thus established the
following continuity result:
\begin{proposition}
\label{pr:continuity}
Let $G,G_1,G_2,\ldots$ be elements of $\cG^*$ whose adjacency
operators are self-adjoint. Let $\mu_{G},\mu_{G_1},\ldots$ denote
the associated spectral measures at their root, and
$m_{G},m_{G_1},\ldots$ the corresponding Cauchy--Stieltjes transforms.
If $\displaystyle G_n \xrightarrow{n\to\infty}^{\cG^*} G$, then
\[
m_{G_n } \xrightarrow{n\to\infty}^{\cH}m_{G } \quad\mbox{and}\quad
\mu_{G_n }\xrightarrow{n\to\infty}^{\cP(\R)}\mu_{G}.
\]
\end{proposition}

As a consequence, when $G,G_1,G_2,\ldots$ are random
elements of $\cG^*$, the same implication holds with all convergences being
replaced by their distributional versions. More precisely, if the law of
$G_n$ converges weakly to that of $G$, then
\[
m_{G_n } \xrightarrow{n\to\infty}^{\cP( \cH) } m_{G }
\quad\mbox{and}\quad \mu_{G_n }\xrightarrow{n\to\infty}^{\cP( \cP
(\R)) } \mu_{G}.
\]

%s3.3 ###
\subsection{Connection with the empirical spectral measure of a finite graph}

In the case of a finite (nonrooted) graph $G_n=(V_n,E_n)$ on $n$
vertices, the adjacency
operator $A_n$ is a particularly simple object: it is bounded and
self-adjoint, and it has exactly $n$ eigenvalues
$\lambda_1(A_n)\geq\cdots\geq\lambda_n(A_n)$ (with
multiplicities), all
of them being real. Moreover, $\ell^2(V_n)\equiv\C^n$ admits an
orthonormal basis of eigenvectors $(b_1,\ldots, b_n)$, a priori
different from the canonical orthonormal basis
$(e_v)_{v\in V_n}$, such that
\[
\forall x\in\C^n\qquad A_n x=\sum_{i=1}^n \lambda_i(A_n)\langle
x,b_i\rangle b_i.
\]
If $(G_n,v)$ denotes the graph $G_n$ when rooted at $v$, the spectral
measure at the root is simply
\[
\mu_{(G_n,v)}=\sum_{i=1}^n |\langle b_i,e_v\rangle|^2 \delta
_{\lambda_i(A_n)}.
\]
In fact $\mu_{(G_n,v)}$ can be interpreted as the local contribution
of vertex $v$ to
the empirical spectral measure $\mu_n$ of $G_n$. Indeed, the above
formula implies
%
%e16 ###
%
\begin{equation}
\label{eq:decompmu}\frac{1}{n}\sum_{v\in
V_n} \mu_{(G_n,v)}=\frac{1}{n}\sum_{i=1}^n\delta_{\lambda_i(A_n)}
= \mu_n.
\end{equation}
Note that the left-hand side can be reinterpreted as the expectation of
$\mu_{(G_n,\OO)}$ under a uniform choice of the root ${\OO}$. More generally,
if $G_n$ is a random graph on $n$ vertices, we
denote by $U(G_n)$ the random element of $\cG^*$ obtained by rooting
$G_n$ at a uniformly chosen vertex,
independently of the random edge-structure. Similarly, we define
$U_2(G_n)$ as
the random element $( (G_n,{\OO}_1),(G_n,{\OO}_2) ) $ in $\cG^* \times
\cG^*$,
where $({\OO}_1,{\OO}_2)$ is a uniformly chosen pair of vertices.
Finally we let $\mu_n$
denote the (random) empirical spectral measure of the adjacency matrix
of $G_n$. With this
notation, we have the following corollary.
\begin{corollary} \label{cor:ESD}
If $U(G_n)$ converges weakly to a rooted GWT $T$ whose degree
distribution $F_*$ has a finite second moment, then
\[
\lim_{n\to\infty} \EE\mu_n = \EE\mu_{T},
\]
where
$\mu_{T}$ denotes the local spectral measure at the root of $T$.
If moreover $U_2(G_n)$ converges weakly to $(T_1,T_2)$, two independent
copies of $T$, then in probability,
\[
\lim_{n\to\infty} \mu_n = \EE\mu_{T}.
\]
\end{corollary}

In the above-mentioned cases of Erd\H{o}s--R\'enyi random graphs and random
graphs with asymptotic degree distribution $F_*$, the assumption on
$U_2(G_n)$ is easily checked. This corollary implies that the
study of the limiting spectral measure of random tree-like graphs
boils down to the study of the local spectral measure at the root of the
limiting GWT. As we have seen, the latter is fully characterized by a
simple RDE
involving its Cauchy--Stieltjes transform. Note, however, that this
result does not give the full statement of Theorem
\ref{th1}(i); the almost sure convergence will be considered later.
\begin{pf*}{Proof of Corollary \ref{cor:ESD}}
By (\ref{eq:decompmu}), we may write for any bounded continuous
function $f\dvtx\R\to\R$,
\[
\EE\int_\R f \,d\mu_n = \frac{1}{n}\sum_{\OO\in V_n}\EE\int_\R f
\,d\mu_{(G_n,\OO)}\xrightarrow{n\to\infty}^{}\EE
\int_\R f \,d\mu_{T},
\]
where the convergence follows from the weak convergence $U(G_n)\to\cT
$ and the continuity result stated in Proposition \ref{pr:continuity}.
This is exactly saying that $\EE\mu_n \to\EE\mu_{T}$.
If, moreover, $U_2(G_n)$ converges weakly to $(T_1,T_2)$, then by the
same argument,
\begin{eqnarray*}
\EE\biggl(\int_\R f \,d\mu_n\biggr)^2 &=& \frac{1}{n^2} \sum_{\OO
_1\in V_n,\OO_2\in V_n}
\EE\biggl(\int_\R f \,d\mu_{(G_n,\OO_1)}\int_\R f \,d\mu_{(G_n,\OO
_2)}\biggr)\\
&\displaystyle \xrightarrow{n\to\infty}^{}&
\biggl(\EE\int_\R f \,d\mu_{T}\biggr)^2,
\end{eqnarray*}
and therefore, the second moment method suffices to conclude that
\[
\int_\R f \,d\mu_n \xrightarrow{n\to\infty}^{P} \EE\int_\R f
\,d\mu_{T},
\]
which is exactly saying that $\mu_n \to \EE\mu_{T}$ in probability.
\end{pf*}

%s3.4 ###
\subsection{Main result: Convergence of the rank}\label{sec:cvrank}

We are now in position to state the main result of this paper. We
consider a sequence of finite random graphs $G_1,G_2,\ldots$ converging
in distribution (once uniformly rooted) to a GWT whose degree distribution
$F_*$ has a finite second moment. As above, $\varphi_* (x) = \sum_k
F_*( k
) x ^k$ denotes the generating function of $F_*$, and we consider the
function
\begin{eqnarray}
&&M\dvtx
x\in[0,1]\mapsto\varphi'_*(1)x\ox+\varphi_*(1-x)+\varphi_*(1-\ox)-1\nonumber\\
&&\eqntext{\mbox{where }\ox=\varphi'_*(1-x)/\varphi'_*(1).}
\end{eqnarray}
Recall that $M'(x)=\varphi_*''(1-x)(\overline{\overline{x}}-x)$ so
that $M(x)$ is a local extremum if and only if $\overline{\overline{x}}=x$.
\begin{theorem}
\label{th:main}
Assume that $U_2(G_n)$ converges weakly to $(\cT_1,\cT_2)$,
two independent copies of a GWT whose degree distribution
$F_*$ has a finite second moment. If the first local extremum of
$M$ is the global maximum, then in probability,
\[
\lim_{n \to\infty} \frac{1}{n}\rank(A_n) = 1-\max_{x\in[0,1]}M(x).
\]
Moreover, a simple sufficient condition for the assumption on $M$ to
hold is that $\varphi_*''$ is log-concave.
\end{theorem}

If the assumption $U_2(G_n)\to(\cT_1,\cT_2)$ is replaced by the weaker
$U(G_n)\to T$, then we only have convergence of the expected rank.

The log-concavity of $\varphi_*''$ is a sufficient condition for the
first local extremum of
$M$ to be a global maximum. Setting $h\dvtx x\mapsto
\overline{\overline{x}}-x$, we find
\[
\forall x\in(0,1)\qquad h''(x)=\frac{\varphi_*''(1-x
)}{\varphi_*'(1)}\frac{\varphi_*''(1-\overline{x}
)}{\varphi_*'(1)}g(x)
\]
with
\[
g(x)=\frac{\varphi_*''(1-x)\varphi_*'''(1-\overline
{x})}{\varphi_*'(1)\varphi_*''(1-\overline{x})}-\frac{\varphi
_*'''(1-x)}{\varphi_*''(1-x)}.
\]
Now, if $\varphi''_*$ is log-concave, then $x\mapsto\varphi
_*'''(x)/\varphi_*''(x)$ is
nonincreasing on $(0,1)$, and therefore, $g$ is decreasing (as the
difference of a decreasing function and a nondecreasing
one). Consequently, $h''$ can vanish at most once on $(0,1)$, hence $h'$
admits at most two zeros on $[0,1]$, and $h$ at most three. The unique root
$x_c$ of $x=\overline{x}$ is always one of them, and if $x$ is another
one, then
so is $\overline{x}$. Therefore, only two cases are possible:
\begin{itemize}
\item[-] Either $x_c$ is the only zero of $h$; then $h(0)>0$ and
$h(1)<0$, so $M$ is maximum at $x_c$,
\item[-] or $h$ admits exactly three zeros $x_-<x_c<x_+$; in this
case the
decreasing function $g$ has to vanish somewhere in $(0,1)$, so $h''$ is
positive and then negative on $(0,1)$. Consequently,
$h$ is decreasing, then increasing, and then decreasing again. In
other words, $M$ is minimum at $x_c$ and maximum at $x_-,x_+$.
\end{itemize}
In both cases, the first local extremum of $M$ is its global maximum.

The remaining part of this
section is devoted to the proof of Theorem \ref{th:main}. First,
recall that $n ^{-1} \rank(A_n) = 1 - \mu_n (\{0\})$. From Corollary
\ref{cor:ESD}, we have in probability,
\[
\limsup_n \mu_n(\{0\}) \leq\EE\mu_{T}(\{0\}).
\]
In order to prove Theorem \ref{th:main}, it is thus sufficient to
establish that
%
%e17 ###
%
\begin{equation}
\label{eq:liminf}
\liminf_n \EE\mu_n ( \{ 0 \} ) \geq\max_{x\in
[0,1]}M(x).
\end{equation}
To do so, we will use the Karp--Sipser leaf removal algorithm, which was
introduced in \cite{karpsipser} to efficiently build a matching (i.e., a
subset of
pairwise disjoint edges) on a finite graph.

For our purposes, the leaf removal algorithm on a locally finite graph
$G =
(V,E)$ can be described as an iterative procedure that constructs two
nondecreasing sequences
$(\cA_t)_{t\geq0}$ and $(\cB_t)_{t\geq0}$
of subsets of $V$ as follows: we start with
\[
\cA_0 = \{v \in V \dvtx\operatorname{deg} ( v ; G ) = 0 \} \quad
\mbox{and}\quad \cB_0 = \varnothing.
\]
Thus, $\cA_0$ is simply the set of isolated vertices in $G$. Then, at
each step
$t\in\N$, we let $G_t$ be the subgraph of $G$ spanned by the vertex-set
$V_t = V \setminus(\cA_t \cup\cB_t\cup P_t)$, where
$P_0=\varnothing$.
We denote by
\[
L_t = \{v \in V_t \dvtx\operatorname{deg} ( v ; G_t ) = 1\}
\]
the set of its leaves. We also introduce the set of vertices that are
adjacent to those leaves,
\[
W_t = \{v \in V_t
\setminus L_t \dvtx\exists u \in L_t, (u v) \in E\}.
\]
We add to $P_{t}$ the set of pairs of adjacent vertices in
$L_t$,
\[
P_{t+1}=P_{t}\cup\{v\in L_t, \exists u\in L_t, (uv)\in E\}.
\]
Then we set
\[
\cA_{t+1} = \cA_t \cup\{ u \in L_t \dvtx\exists v\in W_t, (uv)\in E \}
\quad\mbox{and}\quad \cB_{t+1} = \cB_t \cup W_t.
\]
In words, for any leaf $u$ of $G_t$ whose (unique) neighbor $v$ is not a
leaf, we add $u$ to $\cA_t$ and $v$ to $\cB_t$.
Then to obtain $G_{t+1}$, all nodes in $\cA_{t+1}\cup\cB_{t+1}\cup
P_{t+1}$ are
removed from $G$ (note that to obtain $G_{t+1}$, all leaves from $G_t$
are removed with their adjacent vertices).
If $L_t$ becomes empty, we
have $(\cA_{t+1},\cB_{t+1}) = (\cA_t,\cB_t)$, and the algorithm
stops.
Finally, in the case where the graph $G$ is finite, we
define
%
%e18 ###
%
\begin{equation}
\label{eq:LR}\LR_t(G)=|\cA_{t}(G)|-|\cB_t(G)|.
\end{equation}
Note that for any finite graph $G$, the sequence $(\LR_t(G))_{t\geq
0}$ is nondecreasing.
Note also that the leaf removal algorithm is well defined for a
(possibly infinite) locally finite graph, but the definition
(\ref{eq:LR}) makes sense only for finite graphs.
The lemma below states a connection between these numbers and the rank of
the adjacency matrix of $G$. It was first observed in \cite{bauerg01}, and
a proof can be found in \cite{cdd}.

Although we will not need it here, let us make for completeness the following
observation, which was the original reason why this algorithm was
introduced for finite graphs:
each time a vertex $v$ is added to $\cB_t$, one may arbitrarily associate
it with one of its neighboring leaves $u_v \in\cA_t$. Similarly, for
every vertex $v$ added to~$P_t$, define $u_v$ as its other neighboring
leaf in $P_t$. The edge-set
$\{(vu_v), v \in\cB_t \cup P_t\} $ is then a
matching of $G$, and it is contained in at least one maximum matching
of $G$. Since the graph is finite, the algorithm stops at a finite
time $t^*$. The subgraph of $G$ spanned by the
vertex-set $V\setminus(\cA_{t^*}\cup\cB_{t^*}\cup P_{t^*})$ is a graph
with minimal degree at least $2$ called the core of the graph.
\begin{lemma}
\label{le:LRdimker}
For any finite graph $G$ with adjacency
matrix $A$, and any $t\in\N$,
\[
\dim\operatorname{ker}( A ) \geq|{\LR_t(G)}|.
\]
\end{lemma}
\begin{pf}
Let $u_1 \in L_0 (G)$ be a leaf of $G$ and $v$ its unique neighboring
vertex. Let $G' = G \setminus\{ u_1 , v\}$ and $A(G')$ the adjacency
matrix of $G'$, we have
\[
\dim\operatorname{ker} A (G) = \dim\operatorname{ker} A (G')
\]
(see \cite{bauerg01}). Now, if $\{u_1, \ldots, u_a \} \subset L_0
(G)$, is the set of leaves adjacent to $v$, then $\{u_2, \ldots, u_a\}
$ are isolated vertices in $G'$. The vectors $e_{u_2}, \ldots,
e_{u_a}$ are thus eigenvectors of the kernel of $A'$. By orthogonal
decomposition, we deduce that
\[
\dim\operatorname{ker} A(G) = a - 1 + \dim\operatorname{ker} \bigl( A (
G \setminus\{v, u_1, \ldots, u_a \})\bigr).
\]
By linearity, we obtain that for any integer $t$,
\begin{eqnarray*}
\dim\operatorname{ker} A (G) &=& |\cA_{t}(G)|-|\cB_t(G)| + \dim
\operatorname{ker} \bigl( A \bigl( G
\setminus(\cA_t \cup\cB_t\cup P_t)\bigr)\bigr)\\
&\geq&|\cA_{t}(G)|-|\cB_t(G)|.
\end{eqnarray*}
\upqed\end{pf}

The lower bound (\ref{eq:liminf}) will now follow from the following
proposition.
\begin{proposition}\label{prop:lr}
Let $T$ be a rooted GWT whose degree distribution $F_*$ has a finite
mean. Then
\[
\lim_{t\to\infty} \PP\bigl(\OO\in\cA_t(T)\bigr)-\PP\bigl(\OO\in
\cB_t(T)\bigr) = M(x_0),
\]
where $x_0 \in[0,1]$ is the location of the first local extremum of $M$.
\end{proposition}
\begin{pf}
The argument is close to that appearing in \cite{karpsipser}, Section 4.
For any vertex $\ibf\neq\OO$, we run the leaf removal algorithm on
$\tT_\ibf$ which is the tree $T_\ibf$ with an additional infinite path
starting from $\ibf$. We first compute the corresponding probabilities
$\alpha_t=\PP(\ibf\in\cA_t(\tT_\ibf))$ and $\beta_t=\PP(\ibf
\in
\cB_t(\tT_\ibf))$. For our purpose, adding the infinite path amounts
to increase artificially the degree of the root by $1$: to be a leaf in
$\tT_\ibf$, the root needs to be isolated in $T_{\ibf}$. By
construction, $\ibf$ is in $\cB_t(\tT_\ibf)$ if and only if one of its
children $\kbf$ is in $\cA_{t}(\tT_\kbf)$. Hence if $N$ denotes the
number of children of~$\ibf$, we
have
\[
\beta_t=\EE[1-(1-\alpha_t)^N] = 1-\varphi(1-\alpha_t),
\]
where $\varphi$ is the generating function of $N$ with distribution $F$
given by (\ref{eq:F}). Similarly, $\ibf$ is in $\cA_t(\tT_\ibf)$
if and only if all its
children $\kbf$ are in $\cB_{t-1}(\tT_\kbf)$, so that
$\alpha_t=\varphi(\beta_{t-1})$.
Hence for all $t\geq1$, we have $\alpha_t
=\varphi( 1-\varphi(1-\alpha_{t-1}))$ and $\alpha_0 = 0$.
Since $x\mapsto\varphi( 1-\varphi(1-x))$ is nondecreasing, $\alpha_t$
converges to $\alpha$, the smallest fixed point of the equation
$x=\varphi( 1-\varphi(1-x))$, and $\beta_t$ converges to $\beta=
1-\varphi(1-\alpha)$. Note that $\varphi(x) = \varphi'_*(x)/\varphi
'_*(1)$, where
$\varphi_*$ is the generating function of $F_*$. Hence, with the
notation of
Section \ref{sec:GWT}, we have $\beta= 1- \overline\alpha$, $\alpha=
\overline{\overline{\alpha}}$. In particular, we get $x_0 = \alpha$.

We now compute $\PP(\OO\in\cA_t(T))-\PP(\OO\in
\cB_t(T))$.
Recall that $D(\OO)$ is the set of neighbors of the root $\OO$.
Here are all the possible cases:
\begin{itemize}
\item[-] if $\forall\ibf\in D(\OO)$, $\ibf\in\cB_{t-1}(\tT_\ibf
)$, then
$\OO\in\cA_{t}(T) $;
\item[-] if there exists $\jbf\in D(\OO)\setminus( \cB_{t-1}(\tT
_\jbf) \cup\cA_t(\tT_\jbf)) $
and $\forall\ibf\in
D(\OO)\setminus\jbf$, $\ibf\in\cB_{t-1}(\tT_\ibf)$, then $\OO
\in
\cA_{t}(T) $;
\item[-] if there exists $\ibf\neq\jbf\in D(\OO)$ such that $\ibf
\in
\cA_t(\tT_\ibf)$ and $\jbf\notin\cB_{t-1}(\tT_\jbf)$, then $\OO
\in\cB_{t}(T)$.
\end{itemize}
In all other cases, $\OO\notin\cA_{t}(T)\cup\cB_t(T)$. In summary,
we have
\begin{eqnarray*}
&&\PP\bigl(\OO\in\cA_t(T)\bigr) \\
&&\qquad= \PP\bigl(\forall\ibf\in D(\OO), \ibf\in
\cB_{t-1}(\tT_\ibf)\bigr)\\
&&\qquad\quad{}+
\PP\bigl(\exists\jbf\in D(\OO)\setminus\bigl( \cB_{t-1}(\tT_\jbf)
\cup\cA_t(\tT_\jbf)\bigr), \forall\ibf\in
D(\OO)\setminus\jbf, \ibf\in\cB_{t-1}(\tT_\ibf)\bigr) \\
&&\qquad= \varphi_*(\beta_{t-1})+(1-\beta_{t-1}- \alpha_t)\varphi
_*'(\beta_{t-1}), \\
&&\PP\bigl(\OO\in\cB_t(T)\bigr)\\
&&\qquad= \PP\bigl(\exists\ibf\neq\jbf\in D(\OO
), \ibf\in\cA_t(\tT_\jbf) , \jbf\notin\cB_{t-1}(\tT_\jbf
)\bigr)\\
&&\qquad = \PP\bigl(\exists\ibf\in D(\OO), \ibf\in\cA_t(\tT_\ibf)
\bigr)\\
&&\qquad\quad{} - \PP\bigl( \exists\ibf\in D(\OO), \ibf\in\cA_t(\tT
_\ibf), \forall\jbf\in
D(\OO)\setminus\ibf, \jbf\in\cB_{t-1}(\tT_\jbf)\bigr) \\
&&\qquad=1-\varphi_*(1-\alpha_{t}) - \alpha_t \varphi'_*(\beta_{t-1}).
\end{eqnarray*}
Hence,
\begin{eqnarray*}
\lim_{t\to\infty}\PP\bigl(\OO\in\cA_t(T)\bigr)-\PP\bigl(\OO\in
\cB_t(T)\bigr) &=&
\varphi_*(\beta)+(1-\beta)\varphi_*'(\beta)+\varphi_* (1- \alpha
)-1\\
&=& M(\alpha) = M(x_0),
\end{eqnarray*}
where we have used the identities: $ \beta= 1 - \overline\alpha$,
$\varphi_*' (x) / \varphi_*' (1) = \overline{1 - x} $ and $\overline
{1- \beta} = \alpha$.
\end{pf}
\begin{pf*}{Proof of Theorem \ref{th:main}}
As already pointed out, it is sufficient to prove (\ref{eq:liminf}).
From Lemma \ref{le:LRdimker}, for any integer $t$,
\[
\EE\mu_n (\{0\}) \geq\frac{1}{n}\EE\LR_t(G_n) = \PP\bigl( \OO\in\cA
_t (G_n) \bigr) - \PP\bigl( \OO\in\cB_t (G_n) \bigr),
\]
where $\OO$ is the uniformly drawn root of $U(G_n)$. Note that the events
$\{\OO\in\cA_t ( G_n) \}$ and $\{\OO\in\cB_t ( G_n) \}$ belong to the
$\sigma$-field generated by $(G_n,\OO)_{t+1}$. Thus the convergence of
$U(G_n)$ implies that for any $t\in\N$,
\[
\lim_{n\to\infty}\PP\bigl( \OO\in\cA_t (G_n) \bigr) - \PP\bigl( \OO\in\cB_t
(G_n) \bigr) = \PP\bigl( \OO\in\cA_t (T) \bigr) - \PP\bigl( \OO\in\cB_t (T) \bigr),
\]
where $T$ is a rooted GWT with degree distribution $F_*$ (this is a
standard application of the objective method \cite{aldste}).
%(see also \cite{matcharxiv} for a more elaborate version in the
%context of matchings). We then use Proposition \ref{prop:lr} to
%conclude the proof of Theorem \ref{th:main}.
\end{pf*}

%s4 ###
\section{Conclusion}\label{rem:match}

As explained in the \hyperref[sec1]{Introduction}, the condition on $M$
in Theorem
\ref{th:main} is restrictive,
and the convergence of the rank when this condition is not met (as in
the example described in the \hyperref[sec1]{Introduction}) is left
open. Without any
condition on the function $M$, our work
gives only the following bounds: assume that $U_2(G_n)$ converges
weakly to $(\cT_1,\cT_2)$, two independent copies of a GWT whose
degree distribution $F_*$ has a finite second moment, then in
probability,
%
%e19 ###
%
\begin{eqnarray}
\label{eq:rankbounds}
1-\max_{x\in[0,1]}M(x)&\leq&\liminf_{n\to\infty}\frac{1}{n}\rank
(A_n)\leq\limsup_{n\to\infty}\frac{1}{n}\rank(A_n)\nonumber\\[-8pt]\\[-8pt]
&\leq&1-M(x_0),\nonumber
\end{eqnarray}
where $x_0$ is the first local extremum of $M$.
For example, if the sequence of graphs converges weakly to a GWT
with degree distribution $F_*$ with $F_*(1)=0$, that is, with no leaf,
then $x_0=0$ and $M(0)=F_*(0)$ so that the upper bound in
(\ref{eq:rankbounds}) is trivial.

Our proof for the upper bound on the rank of $A_n$ relies on the
analysis of the leaf removal algorithm on the graph $G_n$. As
explained above, this algorithm when applied to a finite graph produces
a matching and a subgraph of minimal degree~$2$ called the core.
It turns out that the RDEs (\ref{eq:operator}) and (\ref{eq:rdezero})
also appear in the the analysis of the size of maximal matchings on graphs
\cite{matcharxiv}. In particular, if the size of the core is $o(n)$,
the leaf
removal produces an (almost) maximal matching [with error $o(n)$], and
the bounds in (\ref{eq:rankbounds}) match. If the size of the core is
not negligible, but the bounds in (\ref{eq:rankbounds}) match (as,
e.g., in the case where $\varphi_*''$ is log-concave), our result
shows that the asymptotic size of the kernel of the core is zero. In
\cite{matcharxiv}, it is shown that this case corresponds to the situation
where there is an (almost) perfect matching on the core of the graph.
However, as soon as $M(x_0)\neq\max_{x\in[0,1]}M(x)$, for any maximal
matching, there is a positive fraction of vertices in the core that
are not matched \cite{matcharxiv}. In this latter case, the
convergence of the rank is left open.

\begin{appendix}\label{sec:ap}
\section*{\texorpdfstring{Appendix: Proof of Theorem \lowercase{\protect\ref{th1}}}{Appendix: Proof of Theorem 1}}

\begin{figure}[b]

\includegraphics{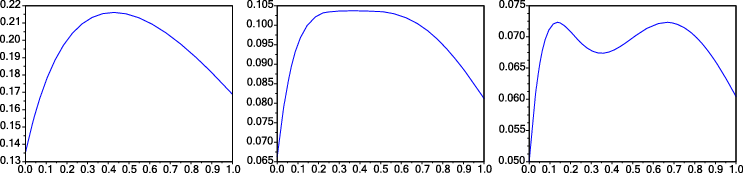}

\caption{From left to right: plot of $M$ for $c = 2$, $c = e$ and $c
= 3$.} \label{fig1}
\end{figure}

In the case where $F_*$ is the Poisson(c) distribution, we simply have
\[
\forall x\in(0,1)\qquad \varphi(x) = \varphi_*(x)=\exp\bigl( c( x - 1) \bigr),
\]
whose second derivative is clearly log-concave. We may therefore apply
Theorem~\ref{th:main} to the sequence of Erd\H{o}s--R\'enyi graphs
$(G_n)_{n \in\N}$.
See Figure \ref{fig1} for a plot of the corresponding function.

%
%f2 ###
%

To complete the proof of Theorem \ref{th1}, it only remains to improve
the convergence in probability into an almost sure convergence. This is
performed by a standard exploration procedure of the edges $E_n$ of the
graph $G_n$. For $1 \leq k \leq n$, we define the random variable in $\{
0,1\}^k$,
\[
X_k = ( A_{i k} )_{ 1 \leq i \leq k }.
\]
By construction, the variables $(X_k)_{1 \leq k \leq n}$ are
independent random variables. Note also that the upper half of the
adjacency matrix $A_n$ is precisely $(X_1, \ldots, X_n)$ and we may
safely write $A_n = A(X_1,\ldots,X_n)$.

For $1 \leq i \leq n$, let $A_{i} (X_1,\ldots,X_n)$ be the principal
minor of $A$ obtained by removing $i$th row and column. If $\lambda
_{1} \leq\cdots\leq\lambda_{n}$ and $\lambda_{1,i} \leq\cdots
\leq\lambda_{n-1,i}$ denote the eigenvalues of $A (X_1,\ldots,X_n)$
and $A_{i} (X_1,\ldots,X_n)$, by the Cauchy interlacing theorem, for
all $1 \leq j \leq n - 1$,
\[
\lambda_j \leq\lambda_{j,i} \leq\lambda_{j+1}.
\]
In particular,
\[
| {\dim\operatorname{ker} A}(X_1,\ldots,X_n) - \dim
\operatorname{ker} A_{i} (X_1,\ldots,X_n) | \leq1.
\]
We note that $A_{i} (X_1,\ldots,X_n)$ does not depend on $X_i$.
Therefore, for all $(x_j \in\{0,1\}^j), 1 \leq j \leq n$, $x'_i \in\{
0,1\}^i$:
\begin{eqnarray*}
&&| {\dim\operatorname{ker} A}(x_1,\ldots,x_{i-1}, x_i , x_{i+1}, \ldots,
x_n)\\
&&\qquad{} - \dim\operatorname{ker} A(x_1,\ldots,x_{i-1}, x'_i ,
x_{i+1}, \ldots, x_n) | \leq2.\nonumber
\end{eqnarray*}
In other words, the function $(x_1,\ldots,x_n) \mapsto\dim
\operatorname{ker} A(x_1,\ldots, x_n)$ is $2$-Lipschitz for the
Hamming distance. By a standard use of Azuma's martingale difference
inequality we get
\[
\PP\bigl( | {\dim\operatorname{ker} A}(X_1,\ldots,X_n) - \EE\dim
\operatorname{ker} A(X_1,\ldots,X_n)| \geq t\bigr) \leq2 \exp
\biggl({-t^2 \over8 n } \biggr).
\]
From the Borel--Cantelli lemma, we obtain that almost surely,
\[
\lim_{n} \frac{\dim\operatorname{ker} A(X_1,\ldots,X_n) - \EE\dim
\operatorname{ker} A(X_1,\ldots,X_n)}{n} = 0.
\]
Since we have already proved that $\EE\dim\operatorname{ker}
A(X_1,\ldots,X_n) / n$ converges to $\max_{ x \in[0,1]} M( x)$, we
deduce that $\dim\operatorname{ker} A(X_1,\ldots,X_n) / n$
convergences a.s. to $\max_{ x \in[0,1]} M( x)$.

It remains to deal with the almost sure convergence in Theorem \ref
{th1}(i). We have already proved that $\mu_n$ converges in probability
to $\mu$. Henceforth $\EE\mu_n$ converges to $\mu$. It is thus
sufficient to prove that almost surely, for all $t \in\R$, $ \mu_n
(( - \infty, t] ) - \EE\mu_n (( - \infty, t] )$ converges to $0$.
The next lemma is a consequence of Lidskii's inequality. For a proof
see Theorem 11.42 in \cite{baisilversteinbook}.
\begin{lemma}[(Rank difference inequality)]\label{le:rankineq}
\label{le:Dn2} Let $A$, $B$ be two $n \times n$ Hermitian matrices
with empirical spectral measures $\mu_A= \frac1 n \sum_{i=1} ^n
\delta_{\lambda_i ( A) }$ and $\mu_B = \frac1 n \sum_{i=1} ^n
\delta_{\lambda_i ( B) }$. Then
\[
\sup_{t \in\R} | \mu_A ((-\infty,t]) - \mu_B ((-\infty,t])
| \leq\frac1 n \rank( A - B).
\]
\end{lemma}

Again, we view $\mu_n$ as a function of $(X_1, \ldots, X_n)$, and
write $\mu_n = \mu_{( X_1, \ldots, X_n)}$. Note that for all $(x_j
\in\{0,1\}^j), 1 \leq j \leq n$, $x'_i \in\{0,1\}^i$, $ A(x_1,\ldots
,x_{i-1}, x_i , x_{i+1}$, $\ldots, x_n) - A(x_1,\ldots,x_{i-1}, x'_i ,
x_{i+1}, \ldots, x_n) $ has only the $i$th row possibly different from
$0$, and we get
\[
\rank\bigl( A(x_1,\ldots,x_{i-1}, x_i , x_{i+1}, \ldots, x_n) -
A(x_1,\ldots,x_{i-1}, x'_i , x_{i+1}, \ldots, x_n) \bigr) \leq2.
\]
Therefore from Lemma \ref{le:Dn2}, for any real $t$,
\[
\bigl| \mu_{(x_1,\ldots,x_{i-1}, x_i , x_{i+1}, \ldots, x_n)}
((-\infty, t]) - \mu_{(x_1,\ldots,x_{i-1}, x'_i , x_{i+1}, \ldots,
x_n)} ((-\infty, t]) \bigr| \leq\frac1 n .
\]
Again, Azuma's martingale difference inequality leads to
\[
\PP\bigl( \bigl| \mu_{(X_1,\ldots, X_n)} ((-\infty, t]) - \EE\mu
_{(X_1,\ldots, X_n)} ((-\infty, t]) \bigr| \geq s \bigr) \leq2
\exp\biggl({- n s^2 \over2 } \biggr).
\]
We deduce similarly from the Borel--Cantelli lemma that $ \mu_n (( -
\infty, t] ) - \EE\mu_n (( - \infty, t] )$ converges a.s. to $0$
and the proof of Theorem \ref{th1} is complete.
\end{appendix}

% imsref loaded by lrinkeviciute, 2010-09-20 10:08:30
%

%
\printaddresses

\end{document}